\documentclass[reqno]{amsart}
\title[On a formula of the $q$-series $_{2k+4}\phi_{2k+3}$]{On a formula of the $q$-series $_{2k+4}\phi_{2k+3}$ and its applications}
\usepackage{amssymb,amsmath,amsthm,epsfig,graphics,latexsym,hyperref}
\usepackage{cite}
\theoremstyle{definition}
\newtheorem{definition}{Definition}
\theoremstyle{plain}

\newtheorem{theorem}    {Theorem}

\newtheorem{corollary}  {Corollary}

\theoremstyle{remark}

\numberwithin{equation}{section}

\newcommand{\fr}{\frac}

\newmuskip\pFqskip
\mathchardef\pFcomma=\mathcode`, 
\newcommand*\pFq[5]{%
  \begingroup
  \begingroup\lccode`~=`,
  \lowercase{\endgroup\def~}{\pFcomma\mkern\pFqskip}%
 \mathcode`,=\string"8000
 {}_{#1}\phi_{#2}\biggl[\genfrac..{0pt}{}{#3}{#4};#5\biggr]%
 \endgroup
 }
 \newmuskip\pGqskip
\mathchardef\pGcomma=\mathcode`, 
\newcommand*\pGq[5]{%
  \begingroup
  \begingroup\lccode`~=`,
    \lowercase{\endgroup\def~}{\pGcomma\mkern\pGqskip}%
  \mathcode`,=\string"8000
  {}_{#1}\psi_{#2}\biggl[\genfrac..{0pt}{}{#3}{#4};#5\biggr]%
  \endgroup
}
\begin{document}
\author[ G. E. Andrews and M. El Bachraoui]{George E. Andrews and Mohamed El Bachraoui}
\address{The Pennsylvania State University, University Park, Pennsylvania 16802}
\email{andrews@math.psu.edu}
\address{Dept. Math. Sci,
United Arab Emirates University, PO Box 15551, Al-Ain, UAE}
\email{melbachraoui@uaeu.ac.ae}
\begin{abstract}
In this paper we apply a formula of the very-well poised $_{2k+4}\phi_{2k+3}$ to write a $k$-tuple sum of $q$-series
as a linear combination of terms wherein each term is a product of expressions of the form
$\fr{1}{(qy, qy^{-1};q)_\infty}$. As an application,
we shall express a variety of sums and double sums of $q$-series as linear combinations of infinite products.
Our formulas are motivated by their connection to overpartition pairs.
\end{abstract}
\date{\textit{\today}}
\thanks{First author partially supported by Simons Foundation Grant 633284}
\maketitle
\section{Introduction}\label{sec introduction}
Throughout $q$ will denote a complex number satisfying $|q|<1$.
Let $\mathbb{N}$ denote the set of positive integers and let
$\mathbb{N}_0=\mathbb{N}\cup\{0\}$ denote the set of nonnegative integers, and let $m, n\in\mathbb{N}_0$.
We adopt the following standard notation from the theory of $q$-series~\cite{Andrews, Gasper-Rahman}
\[
(a;q)_0 = 1,\  (a;q)_n = \prod_{j=0}^{n-1} (1-aq^j),\quad
(a;q)_{\infty} = \prod_{j=0}^{\infty} (1-aq^j),
\]
\[
(a_1,\ldots,a_k;q)_n = \prod_{j=1}^k (a_j;q)_n,\ \text{and\ }
(a_1,\ldots,a_k;q)_{\infty} = \prod_{j=1}^k (a_j;q)_{\infty}.
\]
Furthermore, we will assume that $\fr{1}{(q;q)_{k}} = 0$ for every negative integer $k$.

For complex numbers $y$ and $z$ such that $y\not= z$, let
\[
C(z,y):= \fr{1}{z+z^{-1}-y-y^{-1}}\  \text{and\ }
D(z,y)=(1-y)(1-y^{-1}) (1-z)(1-z^{-1}).
\]
Our primary goal in this paper is to prove the following new result which
enables one to write $k$-tuple sums as linear combinations of terms which are products of expressions
of the form $\fr{1}{(qy, qy^{-1};q)_\infty}$.
\begin{theorem}\label{thm general}
Let $k\in\mathbb{N}$ and $N\in\mathbb{N}_0$. Then
there holds
\[
\sum_{m_1,\ldots,m_{k-1}=0}^\infty
\fr{(b_2,b_2^{-1};q)_{m_1} (b_3,b_3^{-1};q)_{m_1+m_2}\cdots (b_k,b_k^{-1};q)_{m_1+\cdots m_{k-1}} q^{\sum_{i=1}^{k-1}(k-i)m_i}}
{(qb_1,qb_1^{-1};q)_{m_1} (qb_2,qb_2^{-1};q)_{m_1+m_2}\cdots (qb_{k-1},qb_{k-1}^{-1};q)_{m_1+\cdots m_{k-1}}}
\]
\[
 = (qb_k,q/b_k;q)_\infty \prod_{i=1}^k(1-b_i)(1-b_i^{-1}) \sum_{i=1}^k \fr{A_{k,i}(b_1,\ldots,b_k)}{(b_i,b_i^{-1};q)_\infty},
\]
where $A_{1,1}(b_1)=1$ and
\[
A_{k,i}(b_1,\ldots,b_k)
=\begin{cases}
C(b_k,b_{k-1}) A_{k-1,k-1}(b_1,\ldots,b_{k-2},b_{k}) & \text{if $i=k$,} \\
-C(b_k,b_{k-1}) A_{k-1,k-1}(b_1,\ldots,b_{k-1}) & \text{if $i=k-1$,} \\
C(b_k,b_{k-1}) \Big( A_{k-1,i}(b_1,\ldots,b_{k-2},b_{k}) & \\
\qquad\qquad-A_{k-1,i}(b_1,\ldots,b_{k-1})\Big) & \text{if $1\leq i\leq k-2$.}
\end{cases}
\]
\end{theorem}
%
We state the special cases $k=2$ and $k=3$ of Theorem~\ref{thm general} in the following corollary.
\begin{corollary}\label{cor k-2,3}
We have
\[
\begin{split}
{\rm (a)\ }&
\sum_{n\geq 0} \fr{q^n (z,z^{-1};q)_n }{(qy,q/y;q)_n}
= C(z,y)\Big( (1-y)(1-y^{-1}) - \fr{(z,z^{-1};q)_\infty}{(qy,q/y;q)_\infty}\Big). \\
{\rm (b)\ } &
\sum_{m,n\geq 0}\fr{q^{2m+n} (y,y^{-1};q)_m (z,z^{-1};q)_{m+n}}{(qx,q/x;q)_m (qy,q/y;q)_{m+n}}  \\
& =D(x,y)C(z,y)C(z,x)
-D(x,z)C(z,y)C(y,x) \fr{(qz,q/z;q)_\infty}{(qy,q/y;q)_\infty} \\
&+D(y,z)C(z,y)\big( C(y,x)- C(z,x) \big) \fr{(qz,q/z;q)_\infty}{(qx,q/x;q)_\infty}.
\end{split}
\]
\end{corollary}
Our second goal in this paper is to apply Corollary~\ref{cor k-2,3} to express a variety of sums and double sums
of $q$-series as linear combinations of quotients of infinite products.
Note that Corollary~\ref{cor k-2,3}(a) is an equivalent form of Chan and Mao's equation~\cite[(1.4)]{Chan-Mao}.
\begin{definition}
An overpartition~\cite{Corteel-Lovejoy} of $n$ is a partition of $n$ where the first
occurrence of each part may be overlined. The number of overpartitions of $n$, written
$\overline{p}(n)$, has the following generating function
\[
\sum_{n=0}^\infty \overline{p}(n) q^n = \fr{(-q;q)_\infty}{(q;q)_\infty}.
\]
Note that overlined parts in overpartitions are distinct by definition.
We will say that an overpartition has distinct parts if its non-overlined parts are also distinct and we let
$\overline{p}_d(n)$ count the number of such overpartitions. Clearly we have
\[
\sum_{n=0}^\infty \overline{p}_d(n) q^n = (-q;q)_\infty(-q;q)_\infty = (-q;q)_\infty^2.
\]
\end{definition}
\begin{definition}
An overpartition pair~\cite{Lovejoy 2006} of $n$ is a pair of overpartitions $\pi=(\lambda_1, \lambda_2)$ where the sum
of all the parts is $n$. The number of all overpartition pairs of $n$, written $\overline{pp}(n)$, has the following
generating function
\[
\sum_{n=0}^\infty \overline{pp}(n) q^n = \fr{(-q;q)_\infty^2}{(q;q)_\infty^2}.
\]
We say that an overpartition pair $(\lambda_1, \lambda_2)$ has distinct parts if both $\lambda_1$
and $\lambda_2$ have distinct parts and we let $\overline{pp}_d(n)$ denote the number of such overpartition pairs. We clearly have
\[
\sum_{n=0}^\infty \overline{pp}_d(n) q^n = \fr{(-q;q)_\infty^2}{(-q;q)_\infty^2} = (-q;q)_\infty^4.
\]
\end{definition}
In this regard, it is worth to mention that the case $y=-1$ of Corollary~\ref{cor k-2,3}(a) appears in
the paper by Bringmann and Lovejoy~\cite{Bringmann-Lovejoy} as a generating function for
the number $\overline{NN}(n,m)$ of overpartition pairs of $n$ with rank $m$.
It turns out that our $q$-series are related to differences of subclasses of overpartition pairs.
For more on differences of partitions, overpartitions, and overpartitions pairs,
we refer to~\cite{Andrews 2013, Andrews-Yee, Bachraoui 2023-b,  Berkovich-Uncu 2019, Kim-Kim-Lovejoy 2021, Lovejoy 2005, Andrews-Bachraoui}.
We now introduce two examples of overpartition pairs that are related to our $q$-series.
\begin{definition}\label{def A}
For any positive integer $n$ let $A(n)$ denote the number of overpartition pairs $\pi=(\lambda_1, \lambda_2)$ of
$n$ into distinct parts where the smallest part $s(\lambda_1)$ of $\lambda_1$ occurs overlined,
the overlined parts of $\lambda_2$ are $>s(\lambda_1)$ and its non-overlined parts are all $\equiv 0 \bmod{3}$ and $<3s(\lambda_1)$.
Let $A_0(n)$ (resp. $A_1(n)$) denote the number of overpartition pairs counted by $A(n)$ in which the number of
non-overlined parts is even (resp. odd) and let
\[
A'(n) = A_0(n)-A_1(n).
\]
By letting the term $q^{n} (-q^{n+1};q)_{\infty}$ generate the overlined parts of $\lambda_1$
and $(q^{n};q)_{\infty} $ generate its non-overlined parts with the alternating weight $(-1)^k$ of the number of parts and letting
the term $(-q^{n+1};q)_{\infty} $ generate the overlined parts of $\lambda_2$
and $(q^{3};q^3)_{n-1} $ generate its non-overlined parts with the alternating weight of the number of parts,
it is readily verified that
\begin{equation}\label{gen A'}
\sum_{n=1}^\infty A'(n) q^n
=\sum_{n=1}^\infty q^{n} (q^n;q)_\infty (-q^{n+1};q)_{\infty}^2(q^{3};q^3)_{n-1}.
\end{equation}
Furthermore, let $A_2(n)$ (resp. $A_3(n)$) denote the number of overpartition pairs counted by $A(n)$ in which the number of
parts is even (resp. odd) and let
\[
A''(n) = A_3(n)-A_2(n).
\]
It is easy to check that
\begin{equation}\label{gen A''}
\sum_{n=1}^\infty A''(n) q^n
=\sum_{n=1}^\infty q^{n} (q^n;q)_\infty (q^{n+1};q)_{\infty}^2(q^{3};q^3)_{n-1}.
\end{equation}
\end{definition}
For example, $A(5)=14$ counting
\[
\big((\bar{5}),\emptyset \big), \big((\bar{4},\bar{1}),\emptyset \big), \big( (4,\bar{1}),\emptyset \big), \big( (\bar{1}),(\bar{4})\big),
\big((\bar{3},\bar{2}),\emptyset \big), \big((3,\bar{2}),\emptyset \big),
\big( (\bar{2}),(\bar{3}) \big),\big( (\bar{2}),(3) \big),
\]
\[
\big( (\bar{3},\bar{1},1),\emptyset \big), \big( (3,\bar{1},1),\emptyset \big),
\big((\bar{1},1),(\bar{3})\big), \big( (\bar{2},2,\bar{1}), \emptyset \big), \big( (\bar{2},\bar{1}),(\bar{2})\big),
\big( (\bar{2},\bar{1}),(2)\big)
\]
We have $A_0(5)=7$ counting
\[
\big((\bar{5}),\emptyset \big), \big((\bar{4},\bar{1}),\emptyset \big), \big( (\bar{1}),(\bar{4})\big),
\big((\bar{3},\bar{2}),\emptyset \big), \big( (\bar{2}),(\bar{3}) \big), \big( (3,\bar{1},1),\emptyset \big),
\big( (\bar{2},\bar{1}),(\bar{2})\big)
\]
and $A_1(5)=7$ counting the remaining pairs and thus $A'(5)=7-7=0$.
Furthermore, we have $A''(5)=0$ as $A_2(5)=7$ enumerating
\[
\big((\bar{4},\bar{1}),\emptyset \big), \big( (4,\bar{1}),\emptyset \big), \big( (\bar{1}),(\bar{4})\big),
\big((\bar{3},\bar{2}),\emptyset \big), \big((3,\bar{2}),\emptyset \big),
\big( (\bar{2}),(\bar{3}) \big),\big( (\bar{2}),(3) \big).
\]
Our first two identities~\eqref{A'-id} and~\eqref{A''-id} obtained below in Theorem~\ref{thm A'-1} refer to the generating functions of $A'(n)$ and $A''(n)$
and the remaining identities are all companions of these two generating functions. Due to the presence of $(q^3;q^3)_n$ in~\eqref{gen A'}-\eqref{gen A''},
we will need in our proofs the basic fact~\cite[(I.29)]{Gasper-Rahman} that if $\omega=e^{\fr{2\pi i}{3}}$ denotes the third root of unity, then
\[
(x^3;q^3)_n = (x, x\omega, x\omega^{-1};q)_n.
\]

The paper is organized as follows.
In Section~\ref{sec results} we state our main theorems and in
Sections~\ref{sec proof A'-1}-\ref{sec proof B'} we give their proofs.
The proof of Theorem~\ref{thm general} is presented in Section~\ref{sec proof general}.
Finally Section~\ref{sec conclusion} is devoted to some comments and open questions suggested by this work.
\section{The results}\label{sec results}
\begin{theorem}\label{thm A'-1}
We have
\begin{equation}\label{A'-id}
\sum_{n\geq 1} q^n (q^n;q)_\infty (-q^{n+1};q)_{\infty}^2(q^{3};q^3)_{n-1}
= \fr{(q^2;q^2)}{(q;q^2)_\infty} - (q^3;q^3)_\infty.
\end{equation}
\begin{equation}\label{A''-id}
\sum_{n\geq 1} q^n (q^n;q)_\infty (q^{n+1};q)_{\infty}^2(q^{3};q^3)_{n-1}
= \fr{1}{3} \big( (q^3;q^3)_\infty -(q;q)_\infty^3 \big).
\end{equation}
\begin{equation}\label{Comp-1}
\sum_{n\geq 1}q^n (-q^n;q)_\infty (-q^{n+1};q)_{\infty}^2(-q^{3};q^3)_{n-1}
= \fr{1}{3} \big( (-q;q)_\infty^3-(-q^3;q^3)_\infty \big).
\end{equation}
\begin{equation}\label{Comp-2}
\sum_{n\geq 1}q^n (-q^n;q)_\infty (q^{n+1};q)_{\infty}^2(-q^{3};q^3)_{n-1}
=  (-q^3;q^3)_\infty -(-q;q)_\infty(q;q)_\infty^2.
\end{equation}
\end{theorem}
We now collect some other category of companions.
\begin{theorem}\label{thm A'-2}
We have
\[
\begin{split}
{\rm(a)\qquad } &
\sum_{n\geq 1} \fr{q^n}{(q^{n+1};q)_\infty (-q^n;q)_\infty^2 (q^3;q^3)_n}
= \fr{1}{ (q^3;q^3)_\infty}-\fr{(q;q^2)_\infty}{(q^2;q^2)_\infty}. \\
{\rm(b)\qquad } &
\sum_{n\geq 1} \fr{q^n}{(-q^{n+1};q)_\infty (-q^n;q)_\infty^2 (-q^3;q^3)_n}
= \fr{1}{3} \Big( \fr{1}{(-q^3;q^3)_\infty}-\fr{1}{(-q;q)_\infty} \Big). \\
{\rm(c)\qquad  }&
\sum_{n\geq 1} q^n \fr{(q^{n};q)_\infty (q^3;q^3)_{n-1}}{(-q^{n+1};q)_\infty (-q^3;q^3)_n}
= \fr{1}{2}\Big( \fr{(q^3;q^3)_\infty}{(-q^3;q^3)_\infty}-\fr{(q;q)_\infty}{(-q;q)_\infty} \Big). \\
{\rm(d)\qquad } &
\sum_{n\geq 1} q^n \fr{(-q^{n};q)_\infty (-q^3;q^3)_{n-1}}{(q^{n+1};q)_\infty (q^3;q^3)_n}
= \fr{-1}{2}\Big( \fr{(-q^3;q^3)_\infty}{(q^3;q^3)_\infty}-\fr{(-q;q)_\infty}{(q;q)_\infty} \Big).
\end{split}
\]
\end{theorem}
We now focus on double sums. As we will see below in the proofs, our first theorem on double sums
is an instance of Corollary~\ref{cor k-2,3}(b) with
\(
(x,y,z) = (\pm 1,\omega,-\omega)\ \text{or\ } (x,y,z) = (\pm 1,-\omega,\omega),
\)
the second is an instance of Corollary~\ref{cor k-2,3}(b) applied to
\(
(x,y,z) = (\pm \omega,-1,\mp \omega)\ \text{or\ } (x,y,z) = (\pm \omega,\mp \omega,-1),
\)
and the third one is an instance of Corollary~\ref{cor k-2,3}(b) with
\(
(x,y,z) = (1,-1,\pm \omega)\ \text{or\ } (x,y,z) = (1,\pm \omega,-1).
\)
\begin{theorem}\label{thm DS-1}
There holds
\[
\begin{split}
{\rm(a)\qquad} &
\sum_{\substack{m\geq 1\\n\geq 0}}\fr{q^{2m+n} (q^3;q^3)_{m-1} (q;q)_{m+n}(-q^3;q^3)_{m+n-1}}
{(q;q)_{m-1} (q;q)_{m}^2(-q;q)_{m+n-1} (q^3;q^3)_{m+n}} \\
&=\fr{1}{3} \fr{(-q^3;q^3)_\infty}{(-q;q)_\infty(q;q)_\infty^2}
+\fr{1}{6} \fr{(-q^3;q^3)_\infty (q;q)_\infty}{(q^3;q^3)_\infty (-q;q)_\infty}-\fr{1}{2}. \\
{\rm(b)\qquad}&
\sum_{\substack{m\geq 1\\n\geq 0}}\fr{q^{2m+n} (-q^3;q^3)_{m-1} (-q;q)_{m+n}(q^3;q^3)_{m+n-1}}
{(-q;q)_{m-1} (q;q)_{m}^2 (q;q)_{m+n-1} (-q^3;q^3)_{m+n}} \\
&=\fr{1}{3} \fr{(q^3;q^3)_\infty}{(q;q)_\infty^3}
-\fr{1}{2} \fr{(q^3;q^3)_\infty (-q;q)_\infty}{(-q^3;q^3)_\infty (q;q)_\infty}+\fr{1}{6}. \\
{\rm(c)\qquad}&
\sum_{\substack{m\geq 1\\n\geq 0}}\fr{q^{2m+n}(q^3;q^3)_{m-1} (q;q)_{m+n}(-q^3;q^3)_{m+n-1}}
{(q;q)_{m-1} (-q;q)_{m}^2(-q;q)_{m+n-1} (q^3;q^3)_{m+n}} \\
&=\fr{1}{3} \fr{(-q^3;q^3)_\infty}{(-q;q)_\infty^3}
-\fr{1}{2} \fr{(-q^3;q^3)_\infty (q;q)_\infty}{(q^3;q^3)_\infty (-q;q)_\infty}+\fr{1}{6}. \\
{\rm(d)\qquad}&
\sum_{\substack{m\geq 1\\n\geq 0}}\fr{q^{2m+n}(-q^3;q^3)_{m-1} (-q;q)_{m+n}(q^3;q^3)_{m+n-1}}
{(-q;q)_{m-1} (-q;q)_{m}^2 (q;q)_{m+n-1} (-q^3;q^3)_{m+n}} \\
&=\fr{1}{3} \fr{(q^3;q^3)_\infty}{(-q;q)_\infty(q^2;q^2)_\infty}
+\fr{1}{6} \fr{(q^3;q^3)_\infty (-q;q)_\infty}{(-q^3;q^3)_\infty (q;q)_\infty}-\fr{1}{2}.
\end{split}
\]
\end{theorem}
\begin{theorem}\label{thm DS-2}
There holds
\[
\begin{split}
{\rm(a)\qquad}&
\sum_{\substack{m\geq 1\\n\geq 0}}\fr{q^{2m+n}(-q;q)_m (-q;q)_{m-1}^2 (q^3;q^3)_{m+n-1}}
{(-q^3;q^3)_m (q;q)_{m+n-1} (-q;q)_{m+n}^2} \\
&=\fr{1}{6} \fr{(-q;q)_\infty (q^3;q^3)_\infty}{(q;q)_\infty (-q^3;q^3)_\infty}
+\fr{1}{3} \fr{(q^3;q^3)_\infty}{(-q;q)_\infty (q^2;q^2)_\infty}-\fr{1}{2}. \\
{\rm(b)\qquad}&
\sum_{\substack{m\geq 1\\n\geq 0}}\fr{q^{2m+n}(q;q)_m (-q;q)_{m-1}^2 (-q^3;q^3)_{m+n-1}}
{(q^3;q^3)_m (-q;q)_{m+n-1} (-q;q)_{m+n}^2} \\
&=\fr{-1}{2} \fr{(q;q)_\infty (-q^3;q^3)_\infty}{(-q;q)_\infty (q^3;q^3)_\infty}
+\fr{1}{3} \fr{(-q^3;q^3)_\infty}{(-q;q)_\infty^3}+\fr{1}{6}. \\
{\rm(c)\qquad}&
\sum_{\substack{m\geq 1\\n\geq 0}}\fr{q^{2m+n}(q^3;q^3)_{m-1} (-q;q)_{m}(-q;q)_{m+n-1}^2 (q;q)_{m+n}}
{(q;q)_{m-1} (-q^3;q^3)_{m} (q^3;q^3)_{m+n}} \\
&=\fr{1}{6} \fr{(-q;q)_\infty^3}{(-q^3;q^3)_\infty}
-\fr{1}{2} \fr{(-q;q)_\infty (q^2;q^2)_\infty}{(q^3;q^3)_\infty} +\fr{1}{3}. \\
{\rm(d)\qquad}&
\sum_{\substack{m\geq 1\\n\geq 0}}\fr{q^{2m+n}(-q^3;q^3)_{m-1} (q;q)_{m}(-q;q)_{m+n-1}^2 (-q;q)_{m+n}}
{(-q;q)_{m-1} (q^3;q^3)_{m} (-q^3;q^3)_{m+n}} \\
&=\fr{-1}{2} \fr{(-q;q)_\infty (q^2;q^2)_\infty}{(q^3;q^3)_\infty}
+\fr{1}{6} \fr{(-q;q)_\infty^3}{(-q^3;q^3)_\infty} +\fr{1}{3}.
\end{split}
\]
\end{theorem}
\begin{theorem}\label{thm DS-3}
There holds
\[
\begin{split}
{\rm(a)\qquad}&
\sum_{\substack{m\geq 1\\n\geq 0}}\fr{q^{2m+n}(-q;q)_{m-1}^2 (q^3;q^3)_{m+n-1}}
{(q;q)_{m}^2(-q;q)_{m+n}^2 (q;q)_{m+n-1}} \\
&=\fr{1}{12} \fr{(q^3;q^3)_\infty}{(q;q)_\infty^3}
+\fr{1}{4} \fr{(q^3;q^3)_\infty}{(q;q)_\infty (-q;q)_\infty^2}-\fr{1}{3}. \\
{\rm(b)\qquad}&
\sum_{\substack{m\geq 1\\n\geq 0}}\fr{q^{2m+n}(q^3;q^3)_{m-1} (-q;q)_{m+n-1}^2 (q;q)_{m+n}}
{(q;q)_{m-1}(q;q)_{m}^2 (q^3;q^3)_{m+n}} \\
&=\fr{1}{12} \fr{(-q;q)_\infty^2}{(q;q)_\infty^2}
-\fr{1}{3} \fr{(-q;q)_\infty^2 (q;q)_\infty}{(q^3;q^3)_\infty}+\fr{1}{4}. \\
{\rm(c)\qquad}&
\sum_{\substack{m\geq 1\\n\geq 0}}\fr{q^{2m+n}(-q;q)_{m-1}^2 (-q^3;q^3)_{m+n-1}}
{(q;q)_{m}^2(-q;q)_{m+n}^2 (-q;q)_{m+n-1}} \\
&=\fr{1}{4} \fr{(-q^3;q^3)_\infty}{(-q;q)_\infty (q;q)_\infty^2}
+\fr{1}{12} \fr{(-q^3;q^3)_\infty}{(-q;q)_\infty^3}-\fr{1}{3}. \\
{\rm(d)\qquad}&
\sum_{\substack{m\geq 1\\n\geq 0}}\fr{q^{2m+n}(-q^3;q^3)_{m-1} (-q;q)_{m+n-1}^2 (-q;q)_{m+n}}
{(-q;q)_{m-1}(q;q)_{m}^2 (-q^3;q^3)_{m+n}} \\
&=\fr{1}{4} \fr{(-q;q)_\infty^2}{(q;q)_\infty^2}
-\fr{1}{3} \fr{(-q;q)_\infty^3}{(-q^3;q^3)_\infty}+\fr{1}{12}.
\end{split}
\]
\end{theorem}
Noting that the left hand-side of Theorem~\ref{thm A'-1}(a) with $q$ replaced by $q^2$ becomes
\[
\sum_{n\geq 1} q^{2n}(q^{2n};q^2)_\infty (-q^{2+2n};q^2)_\infty (q^6;q^6)_{n-1},
\]
we see that the following theorem is also a companion of Theorem~\ref{thm A'-1}.
\begin{theorem}\label{thm B'}
We have
\[
\begin{split}
{\rm(a)\ } &
\sum_{n\geq 1} q^{2n-1}(q^{2n+1};q^2)_{\infty}(q^{2n};q^2)_{\infty}(q^{2n+3};q^2)_{\infty}(q^{6};q^6)_{n-1} \\
&=\fr{1}{1+q+q^2}(q^6;q^6)_\infty -\fr{1}{1-q^3}(q^2;q^2)_\infty (q;q^2)_\infty^2 . \\
{\rm(b)\ } &
\sum_{n\geq 1} q^{2n-1}(q^{2n+1};q^2)_{\infty}(-q^{2n};q^2)_{\infty}(q^{2n+3};q^2)_{\infty}(-q^{6};q^6)_{n-1} \\
&=\fr{1}{1-q+q^2} \big((-q^6;q^6)_\infty-(-q^2,q^3,q;q^2)_\infty \big). \\
{\rm(c)\ } &
\sum_{n\geq 1} \fr{q^{2n-1}}{(q^{2n-1};q^2)_{\infty}(q^{2n+1};q^2)_{\infty}(q^{2n+2};q^2)_{\infty}(q^{6};q^6)_{n}} \\
&=\fr{1}{1+q+q^2} \Big( \fr{1-q}{(q;q^2)_\infty^2 (q^2;q^2)_\infty} - \fr{1}{(q^6;q^6)_\infty}\Big). \\
{\rm(d)\ } &
\sum_{n\geq 1} \fr{q^{2n-1}}{(q^{2n-1};q^2)_{\infty}(q^{2n+1};q^2)_{\infty}(-q^{2n+2};q^2)_{\infty}(-q^{6};q^6)_{n}} \\
&=\fr{1}{1-q+q^2} \Big( \fr{1}{(-q^2,q^3,q;q^2)_\infty} - \fr{1}{(-q^6;q^6)_\infty}\Big).
\end{split}
\]
\end{theorem}
We closed this section with the following double sum versions of Theorem~\ref{thm B'}.
\begin{theorem}\label{thm DS-4}
There holds
\[
\begin{split}
{\rm(a)\qquad}&
\sum_{\substack{m\geq 1\\n\geq 0}}\fr{q^{4m+2n}(q,q^{-1};q^2)_{m} (q^6;q^6)_{m+n-1}}
{(q^2;q^2)_{m}^2(q^2;q^2)_{m+n-1} (q^3,q;q^2)_{m+n}} \\
&=\fr{1}{3} \fr{(q^6;q^6)_\infty}{(q^2;q^2)_\infty^3} + \fr{q(q^6;q^6)_\infty}{(1+q^3)(q^2;q^2)_\infty (q^3;q^2)_\infty^2}
-\fr{(1-q)^2}{3(1+q+q^2)}. \\
{\rm(b)\qquad}&
\sum_{\substack{m\geq 1\\n\geq 0}}\fr{q^{4m+2n}(q,q^{-1};q^2)_{m} (-q^6;q^6)_{m+n-1}}
{(q^2;q^2)_{m}^2(-q^2;q^2)_{m+n-1} (q^3,q;q^2)_{m+n}} \\
&=\fr{(-q^6;q^6)_\infty}{(-q^2;q^2)_\infty(q^4;q^4)_\infty} - \fr{q(-q^6;q^6)_\infty}{(1-q+q^2)(-q^2;q^2)_\infty (q^3,q;q^2)_\infty}
-\fr{(1-q)^2}{(1-q+q^2)}. \\
{\rm(c)\qquad}&
\sum_{\substack{m\geq 1\\n\geq 0}}\fr{q^{4m+2n}(-q^6;q^6)_{m-1}(q,q^{-1};q^2)_{m+n}(-q^2;q^2)_{m+n}}
{(-q^2;q^2)_{m-1}(q^2;q^2)_{m}^2 (-q^6;q^6)_{m+n}} \\
&=\fr{(q,q^3;q^2)_\infty}{(q^2;q^2)_\infty^2}+\fr{q(-q^2,q,q^{-1};q^2)_\infty}{(1-q+q^2)(-q^6;q^6)_\infty}
-\fr{q}{1-q+q^2}. \\
{\rm(d)\qquad}&
\sum_{\substack{m\geq 1\\n\geq 0}}\fr{q^{4m+2n}(q^6;q^6)_{m-1}(q,q^{-1};q^2)_{m+n}(q^2;q^2)_{m+n}}
{(q^2;q^2)_{m-1}(q^2;q^2)_{m}^2 (q^6;q^6)_{m+n}} \\
&=\fr{1}{3}\fr{(q^6;q^6)_\infty}{(q^2;q^2)_\infty^3}+\fr{1}{3}{q(q;q)_\infty(q^{-1};q^2)_\infty}{(1+q+q^2)(q^6;q^6)_\infty}
-\fr{q}{1+q+q^2}.
\end{split}
\]
\end{theorem}
\section{Proof of Theorem~\ref{thm A'-1}}\label{sec proof A'-1}
(a)\ To prove this part note that
\[
\sum_{n=1}^\infty q^n(q^n;q)_\infty (-q^{n+1};q)_\infty^2 (q^3;q^3)_{n-1}
=(q;q)_\infty (-q;q)_\infty^2 \sum_{n=1}^\infty \fr{(q,\omega q,\omega^{-1}q;q)_{n-1} q^n}{(q;q)_{n-1} (-q;q)_{n}^2}
\]
\[
=(q;q)_\infty (-q;q)_\infty^2 \sum_{n=1}^\infty \fr{(w,w^{-1};q)_{n} q^n}{(1-\omega)(1-\omega^{-1})(-q;q)_{n}^2}
\]
\[
=\fr{1}{3}(q;q)_\infty (-q;q)_\infty^2\sum_{n=1}^\infty \fr{(\omega,\omega^{-1};q)_{n} q^n}{(-q;q)_{n}^2}.
\]
Then by an appeal to Corollary~\ref{cor k-2,3}(a) with $z=\omega$ and $y=-1$ we derive
\[
\sum_{n=1}^\infty q^n(q^n;q)_\infty (-q^{n+1};q)_\infty^2 (q^3;q^3)_{n-1}
\]
\[
\fr{1}{3}(q;q)_\infty (-q;q)_\infty^2 \Big(4-\fr{(\omega,\omega^{-1};q)_\infty}{(-q,-q;q)_\infty} \Big)
-\fr{1}{3}(q;q)_\infty (-q;q)_\infty^2
\]
\[
=\fr{1}{3}(q;q)_\infty (-q;q)_\infty^2 \Big(4-\fr{3(q^3;q^3)_\infty}{(q;q)_\infty (-q;q)_\infty^2} \Big)
-\fr{1}{3}(q;q)_\infty (-q;q)_\infty^2
\]
\[
=(q,-q,-q;q)_\infty - (q^3;q^3)_\infty,
\]
which is the desired formula.

(b)\
We now apply Corollary~\ref{cor k-2,3}(a) with $(z,y)=(\omega,1)$.
Then
\[
\sum_{n=1}^\infty q^n(q^n;q)_\infty (q^{n+1};q)_\infty^2 (q^3;q^3)_{n-1}
= \sum_{n=1}^\infty q^n(q^n;q)_\infty (q^{n+1};q)_\infty^2 (q^3;q^3)_{n-1}
\]
\[
=(q;q)_\infty^3 \sum_{n=1}^\infty \fr{(q,\omega q,\omega^{-1}q;q)_{n-1} q^n}{(q;q)_{n-1} (q;q)_{n}^2}
\]
\[
=\fr{(q;q)_\infty^3}{3}\sum_{n=1}^\infty \fr{(\omega,\omega^{-1};q)_{n} q^n}{(q;q)_{n}^2}
\]
\[
=\fr{(q;q)_\infty^3}{3} \Big( -1 + \sum_{n=0}^\infty \fr{(\omega,\omega^{-1};q)_{n} q^n}{(q;q)_{n}^2} \Big),
\]
which by Corollary~\ref{cor k-2,3}(a) with $(z,y)=(\omega,1)$ yields
\[
\sum_{n=1}^\infty A''(n)q^n
=\fr{(q;q)_\infty^3}{3} \Big( -1 + \fr{1}{3}\fr{(\omega,\omega^{-1};q)_\infty}{(q,q;q)_\infty} \Big)
\]
\[
=\fr{1}{3} \big( (q^3;q^3)_\infty - (q;q)_\infty^3 \big).
\]
This confirms the identity.

(c)\ We have
\[
\sum_{n=1}^\infty q^n(-q^n;q)_\infty (-q^{n+1};q)_\infty^2 (-q^3;q^3)_{n-1}
\]
\[
=(-q;q)_\infty^3 \sum_{n=1}^\infty \fr{(-q,-\omega q,-q/\omega;q)_{n-1} q^n}{(-q;q)_{n-1} (-q;q)_{n}^2}
\]
\[
=(-q;q)_\infty^3\sum_{n=1}^\infty \fr{(-\omega,-\omega^{-1};q)_{n} q^n}{(-q;q)_{n}^2}
\]
\[
=(-q;q)_\infty^3 \Big( -1 + \sum_{n=0}^\infty \fr{(-\omega,-\omega^{-1};q)_{n} q^n}{(-q;q)_{n}^2} \Big),
\]
\[
=(-q;q)_\infty^3 \Big(-1 + \fr{1}{3}\big( 4-\fr{(-\omega,-\omega^{-1};q)_\infty}{(-q;q)_\infty^2}\big) \Big)
\]
\[
=\fr{1}{3} \big( (-q;q)_\infty^3 - (-q^3;q^3)_\infty \big),
\]
where the penultimate formula follows by Corollary~\ref{cor k-2,3}(a) applied to $(z,y)=(-\omega,-1)$.

(d)\ This part follows similarly from Corollary~\ref{cor k-2,3}(a) applied to $(z,y)=(-\omega,1)$.
\section{Proof of Theorem~\ref{thm A'-2}}\label{sec proof A'-2}
(a)\ By virtue of Corollary~\ref{cor k-2,3}(a) with $(z,y)=(-1,\omega)$ we find
\[
\sum_{n=1}^\infty \fr{q^n}{(q^{n+1};q)_\infty (-q^n;q)_\infty^2 (q^3;q^3)_n}
\]
\[
=\fr{1}{(q;q)_\infty (-q;q)_\infty^2} \sum_{n=1}^\infty \fr{q^n(q;q)_n (-q;q)_{n-1}^2}{(q,\omega q,\omega^{-1}q;q)_n}
\]
\[
=\fr{1}{4 (q;q)_\infty (-q;q)_\infty^2} \sum_{n=1}^\infty q^n \fr{(-1,-1;q)_n}{(\omega q,\omega^{-1}q;q)_n}
\]
\[
=\fr{1}{4 (q;q)_\infty (-q;q)_\infty^2} \Big( -1 -3 + \fr{(-1;q)_\infty^2}{(\omega q, \omega^{-1}q;q)_\infty} \Big)
\]
\[
=-\fr{(q;q^2)_\infty}{(q^2;q^2)_\infty} + \fr{1}{(q^3;q^3)_\infty},
\]
which is the desired formula.

(b)\ The formula in this part follows similarly by using Corollary~\ref{cor k-2,3}(a) with $(z,y)=(-1,-\omega)$

(c)\ By employing Corollary~\ref{cor k-2,3}(a) with $(z,y)=(\omega,-\omega)$, we find
\[
\sum_{n=1}^\infty q^n \fr{(q^{n};q)_\infty (q^3;q^3)_{n-1}}{(-q^{n+1};q)_\infty (-q^3;q^3)_n}
\]
\[
=\fr{(q;q)_\infty}{(-q;q)_\infty} \sum_{n=1}^\infty q^n
\fr{(-q;q)_n (q,\omega q, \omega^{-1}q;q)_{n-1}}{(q;q)_{n-1} (-q,-\omega q,- \omega^{-1}q;q)_{n}}
\]
\[
=\fr{(q;q)_\infty}{3 (-q;q)_\infty} \sum_{n=1}^\infty q^n \fr{(\omega ,\omega^{-1} ;q)_n}{(-\omega q,-\omega^{-1}q;q)_n}
\]
\[
=\fr{(q;q)_\infty}{3 (-q;q)_\infty} \Big( -1 -\fr{1}{2} + \fr{1}{2}\fr{(\omega ,\omega^{-1} ;q)_\infty}{(-\omega q,-\omega^{-1}q;q)_\infty} \Big)
\]
\[
=-\fr{(q;q)_\infty}{2(-q;q)_\infty}+ \fr{(q^3;q^3)_\infty}{2(-q^3;q^3)_\infty}.
\]
This gives the desired identity.

(d)\ Finally, the proof of this part follows similarly by Corollary~\ref{cor k-2,3}(a)
applied to $(z,y)=(-\omega,\omega)$.
\section{Proof of Theorem~\ref{thm DS-1}}\label{sec proof DS-1}
(a)\ We want to apply Corollary~\ref{cor k-2,3}(b) to $(x,y,z)=(1,\omega,-\omega)$.
Note first that by Corollary~\ref{cor k-2,3} we get
\begin{equation}\label{ds1-a1}
\sum_{n\geq 0}\fr{(-\omega,-\omega^{-1};q)_n q^n}{(q\omega,q/\omega;q)_{n}}
=\fr{3}{2}-\fr{1}{2}\fr{(-q^3;q^3)_\infty (q;q)_\infty}{(q^3;q^3)_\infty (-q;q)_\infty}
\end{equation}
and
\begin{equation}\label{ds1-a2}
\sum_{m,n\geq 0}\fr{(\omega,\omega^{-1};q)_m(-\omega,-\omega^{-1};q)_{m+n} q^{2m+n}}{(q;q)_m^2(q\omega,q/\omega;q)_{m+n}}
=\fr{(-q^3;q^3)_\infty}{(-q;q)_\infty (q;q)_\infty^2}.
\end{equation}
Then we get the desired formula by using~\eqref{ds1-a1},~\eqref{ds1-a2}, and the following equality
\[
\sum_{\substack{m\geq 1\\n\geq 0}}\fr{(q^3;q^3)_{m-1} (q;q)_{m+n}(-q^3;q^3)_{m+n-1}}
{(q;q)_{m-1} (q;q)_{m}^2(-q;q)_{m+n-1} (q^3;q^3)_{m+n}}q^{2m+n}
\]
\[
=\fr{1}{3} \sum_{\substack{m\geq 1\\n\geq 0}}
\fr{(\omega,\omega^{-1};q)_m (-\omega,-\omega^{-1};q)_{m+n}}{(q;q)_m^2 (q\omega,q/\omega;q)_{m+n}} q^{2m+n}
\]
\[
=\fr{1}{3} \Big( \sum_{m,n\geq 0}\fr{(\omega,\omega^{-1};q)_m(-\omega,-\omega^{-1};q)_{m+n} q^{2m+n}}{(q;q)_m^2(q\omega,q/\omega;q)_{m+n}}
-\sum_{n\geq 0}\fr{(-\omega,-\omega^{-1};q)_n q^n}{(q\omega,q/\omega;q)_{n}} \Big).
\]
Parts~(b), (c), and (d) follow similarly from Corollary~\ref{cor k-2,3}(b) applied to $(x,y,z)=(1,-\omega,\omega)$, $(x,y,z)=(-1,\omega,-\omega)$
and to $(x,y,z)=(-1,-\omega,\omega)$ respectively.
\section{Proof of Theorem~\ref{thm DS-2} and Theorem~\ref{thm DS-3}}\label{sec proof DS-2-3}
{\emph Proof of Theorem~\ref{thm DS-2}.\ }
(a)\ Note first that by Corollary~\ref{cor k-2,3} we get
\begin{equation}\label{ds2-a1}
\sum_{n\geq 0}\fr{(\omega,\omega^{-1};q)_n q^n}{(-q;q)_{n}^2}
= 4- \fr{(\omega,\omega^{-1};q)_\infty}{(-q;q)_\infty^2}
=4-3\fr{(q^3;q^3)_\infty}{(-q;q)_\infty (q^2;q^2)_\infty}
\end{equation}
and
\[
\sum_{m,n\geq 0}\fr{(-1,-1;q)_m(\omega,\omega^{-1};q)_{m+n} q^{2m+n}}{(-q\omega,-q/\omega;q)_m (-q;q)_{m+n}^2}
= \fr{2(q\omega,q/\omega;q)_\infty}{(-q\omega,-q/\omega;q)_\infty}
+\fr{(q\omega,q/\omega;q)_\infty}{(-q,-q;q)_\infty} -2
\]
\begin{equation}\label{ds2-a2}
=\fr{2(-q;q)_\infty(q^3;q^3)_\infty}{(q;q)_\infty(-q^3;q^3)_\infty}
+\fr{(q^3;q^3)_\infty}{(-q;q)_\infty(q^2;q^2)_\infty} -2.
\end{equation}
Now by Corollary~\ref{cor k-2,3}(b) applied to $(x,y,z)=(-\omega,-1,\omega)$, we find
\[
\sum_{\substack{m\geq 1\\n\geq 0}}\fr{(-q;q)_m (-q;q)_{m-1}^2 (q^3;q^3)_{m+n-1}}
{(-q^3;q^3)_m (q;q)_{m+n-1} (-q;q)_{m+n}^2}q^{2m+n}
\]
\[
=\fr{1}{4(1-\omega)(1-\omega^{-1})}
\sum_{\substack{m\geq 1\\n\geq 0}}\fr{(-1,-1;q)_{m} (\omega,\omega^{-1};q)_{m+n}}
{(-q\omega,-q/\omega;q)_m (-q,-q;q)_{m+n}}q^{2m+n}
\]
\[
\fr{-1}{12}\sum_{n\geq 0}\fr{(\omega,\omega^{-1};q)_n q^n}{(-q;q)_{n}^2}
+\fr{1}{12} \sum_{m,n\geq 0}\fr{(-1,-1;q)_m(\omega,\omega^{-1};q)_{m+n} q^{2m+n}}{(-q\omega,-q/\omega;q)_m (-q;q)_{m+n}^2}
\]
\[
=\fr{1}{6} \fr{(-q;q)_\infty (q^3;q^3)_\infty}{(q;q)_\infty (-q^3;q^3)_\infty}
+\fr{1}{3} \fr{(q^3;q^3)_\infty}{(-q;q)_\infty (q^2;q^2)_\infty}-\fr{1}{2},
\]
where the last identity follows from~\eqref{ds2-a1} and~\eqref{ds2-a2}.

(b)\ By Corollary~\ref{cor k-2,3} we have
\begin{equation}\label{ds2-b1}
\sum_{n\geq 0}\fr{(-\omega,-\omega^{-1};q)_n q^n}{(-q;q)_{n}^2}
=\fr{4}{3}-\fr{1}{3}\fr{(-q^3;q^3)_\infty}{(-q;q)_\infty^3}
\end{equation}
and
\[
\sum_{m,n\geq 0}\fr{(-1,-1;q)_m(-\omega,-\omega^{-1};q)_{m+n} q^{2m+n}}{(q\omega,q/\omega;q)_m (-q;q)_{m+n}^2}
\]
\begin{equation}\label{ds2-b2}
=\fr{-2(-q^3;q^3)_\infty (q;q)_\infty}{(q^3;q^3)_\infty (-q;q)_\infty}
+\fr{(-q^3;q^3)_\infty}{(-q;q)_\infty^3} +2.
\end{equation}
Then the desired formula follows by combining~\eqref{ds2-b1} and~\eqref{ds2-b2} with the following equality
\[
\sum_{\substack{m\geq 1\\n\geq 0}}\fr{(q;q)_m (-q;q)_{m-1}^2 (-q^3;q^3)_{m+n-1}}
{(q^3;q^3)_m (-q;q)_{m+n-1} (-q;q)_{m+n}^2}q^{2m+n}
\]
\[
\fr{1}{4}\sum_{m,n\geq 0}\fr{(-1,-1;q)_m(-\omega,-\omega^{-1};q)_{m+n} q^{2m+n}}{(q\omega,q/\omega;q)_m (-q;q)_{m+n}^2}
-\fr{1}{4}\sum_{n\geq 0}\fr{(-\omega,-\omega^{-1};q)_n q^n}{(-q;q)_{n}^2}.
\]
Parts~(c) and~(d) follow similarly by Corollary~\ref{cor k-2,3} applied to $(x,y,z)=(-\omega,\omega,-1)$, and $(x,y,z)=(\omega,-\omega,-1)$
respectively.

{\emph Proof of Theorem~\ref{thm DS-3}.\ }
The details of the proofs are omitted since the steps are similar to the above proofs.
Part~(a), (b), (c), and (d) follow similarly from Corollary~\ref{cor k-2,3}(b) applied to $(x,y,z)=(1,-1,\omega)$,
$(x,y,z)=(1,\omega,-1)$, $(x,y,z)=(1,-1,-\omega)$, and $(x,y,z)=(1,-\omega,-1)$ respectively.

\section{Proof of Theorem~\ref{thm B'} and Theorem~\ref{thm DS-4}}\label{sec proof B'}
{\emph Proof of Theorem~\ref{thm B'}.\ }
As for (a), we have
\[
\sum_{n=1}^\infty q^{2n}(q^{2n+1};q^2)_{\infty}(q^{2n};q^2)_{\infty}(q^{2n+3};q^2)_{\infty}(q^{6};q^6)_{n-1}
\]
\[
=(q^2;q^2)_\infty(q;q^2)_\infty^2 \sum_{n=1}^\infty \fr{(q^6;q^6)_{n-1} q^{2n}}{(q^2;q^2)_{n-1} (q;q^2)_n  (-q;q^2)_{n+1}}
\]
\[
=\fr{(q^2;q^2)_\infty(q;q^2)_\infty^2}{1-q} \sum_{n=1}^\infty \fr{(q^2 \omega, q^2\omega^{-1};q^2)_{n-1} q^{2n}}{(q;q^2)_n  (q^3;q^2)_{n}}
\]
\[
=\fr{(q^2;q^2)_\infty(q;q^2)_\infty^2}{3(1-q)} \sum_{n=1}^\infty \fr{(\omega, \omega^{-1};q^2)_{n-1} q^{2n}}{(q^3,q;q^2)_n}.
\]
Then by Corollary~\ref{cor k-2,3}(a) applied to $(z,y)=(\omega,q)$ we get
\[
\sum_{n=1}^\infty q^{2n}(q^{2n+1};q^2)_{\infty}(q^{2n};q^2)_{\infty}(q^{2n+3};q^2)_{\infty}(q^{6};q^6)_{n-1}
\]
\[
=\fr{(q^2;q^2)_\infty(q;q^2)_\infty^2}{3(1-q)}
\Big( -1 +\fr{(1+q)^2}{1+q+q^2} - \fr{q}{1+q+q^2}\fr{(\omega,\omega^{-1};q^2)_\infty}{(q,q^3;q^2)_\infty} \Big)
\]
\[
=\fr{-q}{1-q^3}(q^2;q^2)_\infty (q;q^2)_\infty^2 + \fr{q}{1+q+q^2}(q^6;q^6)_\infty,
\]
as desired.
Parts~(b),~(c), and~(d) follow similarly by Corollary~\ref{cor k-2,3}(a) applied to $(z,y)=(-\omega,q)$,
$(z,y)=(q,\omega)$, and $(z,y)=(q,-\omega)$ respectively.

{\emph Proof of Theorem~\ref{thm DS-4}.\ }
Part~(a) follows from Corollary~\ref{cor k-2,3}(b) applied to $(x,y,z)=(1,q,\omega)$, $(x,y,z)=(1,q,-\omega)$,
$(x,y,z)=(1,-\omega,q)$, and $(x,y,z)=(1,\omega,q)$.
\section{Proof of Theorem~\ref{thm general}}\label{sec proof general}
By~\cite[Theorem 4]{Andrews 1975} applied to $c_i=b_i^{-1}$ we get
\[
\pFq{2k+4}{2k+3}{a,q\sqrt{a},-q\sqrt{a},b_1,b_{1}^{-1},\ldots,b_k,b_{k}^{-1},q^{-N}}
{\sqrt{a},-\sqrt{a},aqb_1,aq/b_1,\ldots,aqb_k,aq/b_k,aq^{N+1}}
{q, a^k q^{k+N}}
\]
\begin{equation}\label{Andr-1}
=\fr{(aq)_N^2}{(aqb_k)_N (aq/b_k)_N}
\sum_{m_1,\ldots,m_{k-1}=0}^\infty
\fr{\prod_{i=1}^{k-1} (aq)_{m_i}}{\prod_{i=1}^{k-1} (q)_{m_i}}
\end{equation}
\[
\cdot \fr{(b_2, b_2^{-1})_{m_1} (b_3,b_3^{-1})_{m_1+m_2}\cdots (b_k,b_k^{-1})_{m_1+\cdots m_{k-1}}}
{(aqb_1,aqb_1^{-1})_{m_1} (aqb_2,aqb_2^{-1})_{m_1+m_2} \cdots (aqb_{k-1},aqb_{k-1}^{-1})_{m_1+\cdots m_{k-1}}}
\]
\[
\cdot \fr{(q^{-N})_{m_1+\cdots+m_{k-1}}}{(q^{-N}/a)_{m_1+\cdots+m_{k-1}}}
(aq)^{m_{k-2}+2m_{k-3}+\cdots+(k-2)m_1} q^{m_1+\cdots m_{k-1}}.
\]
It is easy to see that
\begin{equation}\label{rhs}
\lim_{\substack{a\to 1\\N\to\infty}}\text{R.H.S\eqref{Andr-1}}
\end{equation}
\[
=
\fr{(q;q)_\infty^2}{(qb_k,q/b_k;q)_\infty}
\sum_{m_1,\ldots,m_{k-1}=0}^\infty
\fr{(b_2,b_2^{-1})_{m_1} (b_3,b_3^{-1})_{m_1+m_2}\cdots (b_k,b_k^{-1})_{m_1+\cdots m_{k-1}} q^{\sum_{i=1}^{k-1}(k-i)m_i}}
{(qb_1,qb_1^{-1})_{m_1} (qb_2,b_2^{-1})_{m_1+m_2}\cdots (qb_{k-1},qb_{k-1}^{-1})_{m_1+\cdots m_{k-1}}}.
\]
We now focus on the left hand-side of~\eqref{Andr-1} as $a\to 1$ and $N\to\infty$.
For convenience, we define
\[
L_{k,N}(b_1,\ldots,b_k)=\lim_{a\to 1}\text{L.H.S\eqref{Andr-1}},\
L_k(b_1,\ldots,b_k) =\lim_{N\to\infty} L_{k,N}(b_1,\ldots,b_k),
\]
and
\[
F_k(b_1,\ldots,b_k)= \fr{L_k(b_1,\ldots,b_k)}{\prod_{i=1}^k (1-b_i)(1-b_i^{-1})}.
\]
Then
we have
\[
L_{k,N}(b_1,\ldots,b_k)
=
\lim_{a\to 1}
\sum_{n=0}^N \fr{(a,q\sqrt{a},-q\sqrt{a},b,b^{-1},\ldots,b_k,b_k^{-1},q^{-N};q)_n a^k q^{(k+N)n}}
{(q,\sqrt{a},-\sqrt{a},aqb_1,aqb_1^{-1},\ldots,aqb_k,aqb_k^{-1},aq^{N+1};q)_n}
\]
\[
= 1+ 2\sum_{n=1}^N \fr{(-q,b_1,b_1^{-1},\ldots,b_k,b_k^{-1},q^{-N};q)_n q^{(k+N)n}}
{(-1,qb_1,q/b_1,\ldots,qb_k,q/b_k,q^{N+1};q)_n}
\]
\[
=1+
\sum_{n=1}^N (1+q^n) \fr{\prod_{i=1}^k(1-b_i)(1-b_{i}^{-1})}
{\prod_{i=1}^k(1-b_i q^n)(1-b_{i}^{-1}q^n)} \fr{(q^{-N};q)_n}{(q^{N+1};q)_n} q^{(k+N)n}
\]
\begin{equation}\label{N-finite}
=1 + \prod_{i=1}^k(1-b_i)(1-b_{i}^{-1})
\sum_{n=1}^N (1+q^n)\fr{q^{(k+N)n} (q^{-N};q)_n}
{\prod_{i=1}^k(1-b_i q^n)(1-b_{i}^{-1}q^n)(q^{N+1};q)_n}.
\end{equation}
Now, as for $n\geq 1$
\[
\lim_{N\to\infty} (q^{N+1};q)_n = 1\ \text{and\ } \lim_{N\to\infty}q^{Nn}(q^{-N};q)_n = (-1)^n q^{{n\choose 2}},
\]
we get by letting $N\to\infty$ in~\eqref{N-finite},
\[
L_k(b_1,\ldots,b_k) =
1+ \prod_{i=1}^k (1-b_i)(1-b_i^{-1}) \sum_{n=1}^\infty (1+q^n)
\fr{(-1)^n q^{{n+1\choose 2}+(k-1)n}}{\prod_{i=1}^k (1-b_i q^n) (1-b_i^{-1} q^n)}.
\]
Thus,
\[
F_k(b_1,\ldots,b_k)
=\fr{1}{\prod_{i=1}^k (1-b_i)(1-b_i^{-1})}
+ \sum_{n=1}^\infty (1+q^n)\fr{(-1)^n q^{{n+1\choose 2}+(k-1)n}}{\prod_{i=1}^k (1-b_i q^n) (1-b_i^{-1}q^n)}
\]
\begin{equation}\label{help-id-0}
=\sum_{n=-\infty}^\infty \fr{(-1)^n q^{{n+1\choose 2}+(k-1)n}}{\prod_{i=1}^k (1-b_i q^n) (1-b_i^{-1}q^n)}
\end{equation}
\[
=\sum_{n=-\infty}^\infty \fr{(-1)^n q^{{n+1\choose 2}+(k-2)n}}{\prod_{i=1}^{k-2} (1-b_i q^n) (1-b_i^{-1}q^n)}
\]
\[
\cdot C(b_k,b_{k-1})
\Big( \fr{1}{(1-b_k q^n)(1-b_k^{-1} q^n)} - \fr{1}{(1-b_{k-1} q^n)(1-b_{k-1}^{-1} q^n)} \Big)
\]
\begin{equation}\label{recurrence}
=C(b_k,b_{k-1}) \Big( F_{k-1}(b_1,\ldots,b_{k-2},b_k)-F_{k-1}(b_1,\ldots,b_{k-1}) \Big).
\end{equation}
We will now evaluate $F_1(b_1)$. To this end, we need
Bailey's sum of the well-poised $_3\psi_3$~\cite[(II.31)]{Gasper-Rahman}
\begin{equation}\label{Bailey}
\pGq{3}{3}{b,c,d}{q/b,q/c,q/d}{q, \fr{q}{bcd}}
=\fr{(q,q/bc,q/bd,q/cd)_\infty}{(q/b,q/c,q/d,q/bcd)_\infty}.
\end{equation}
From~\eqref{help-id-0}, we find for $k=1$,
\[
(1-b_1)(1-b_1^{-1}) F_1(b_1)
=\sum_{n=-\infty}^\infty \fr{(1-b_1)(1-b_1^{-1})}{(1-b_1 q^n) (1-b_1^{-1}q^n)} (-1)^n q^{{n+1\choose 2}}
\]
\[
=\sum_{n=-\infty}^\infty \fr{(b_1,b_1^{-1};q)_n}{(b_1q,b_1^{-1}q;q)_n} (-1)^n q^{{n+1\choose 2}}
\]
\begin{equation}\label{help-id-1}
=\fr{(q;q)_\infty^2}{(qb_1;q)_\infty (qb_1^{-1};q)_\infty},
\end{equation}
where the last identity follows from~\eqref{Bailey} applied to $b=c^{-1}=b^{-1}$ and $d\to\infty$.
Now use~\eqref{help-id-1} to deduce
\begin{equation}\label{basic}
F_1(b_1) = \fr{(q;q)_\infty^2}{(1-b_1)(1-b_1^{-1})(qb;q)_\infty (qb_1^{-1};q)_\infty}
=\fr{(q;q)_\infty^2}{(b_1,b_1^{-1};q)_\infty},
\end{equation}
showing that $A_{1,1}(b_1)=1$. Furthermore, combining the recurrence relation~\eqref{recurrence} with the basic case~\eqref{basic} yields
\[
F_k(b_1,\ldots,b_k)
= C(b_k,b_{k-1})
\Big( \sum_{i=1}^{k-1} \fr{A_{k-1,i}(b_1,\ldots,b_{k-2},b_k)}{(b_i,b_i^{-1};q)_\infty}
\]
\[
-\sum_{i=1}^{k-1} \fr{A_{k-1,i}(b_1,\ldots,b_{k-1})}{(b_i,b_i^{-1};q)_\infty} \Big)
\]
\[
=C(b_k,b_{k-1})
\Big( \sum_{i=1}^{k-2} \fr{A_{k-1,i}(b_1,\ldots,b_{k-2},b_k)}{(b_i,b_i^{-1};q)_\infty}
\]
\[
+ \fr{A_{k-1,k-1}(b_1,\ldots,b_{k-2},b_k)}{ (b_k,b_k^{-1};q)_\infty}
-\sum_{i=1}^{k-1} \fr{A_{k-1,i}(b_1,\ldots,b_{k-1})}{(b_i,b_i^{-1};q)_\infty} \Big).
\]
This implies that
\[
\begin{split}
A_{k,k}(b_1,\ldots,b_k)
&= C(b_k,b_{k-1})A_{k-1,k-1}(b_1,\ldots,b_{k-2},b_k), \\
A_{k,k-1}(b_1,\ldots,b_k)
&= -C(b_k,b_{k-1})A_{k-1,k-1}(b_1,\ldots,b_{k-1}),
\end{split}
\]
and for $1\leq i\leq k-2$,
\[
A_{k,i}(b_1,\ldots,b_k)
= C(b_k,b_{k-1}) \Big(A_{k-1,i}(b_1,\ldots,b_{k-2},b_k) - A_{k-1,i}(b_1,\ldots,b_{k-1}) \Big).
\]
That is,
\begin{equation}\label{A k-i}
A_{k,i}(b_1,\ldots,b_k)
=\begin{cases}
C(b_k,b_{k-1}) A_{k-1,k-1}(b_1,\ldots,b_{k-2},b_{k}) & \text{if $i=k$,} \\
-C(b_k,b_{k-1}) A_{k-1,k-1}(b_1,\ldots,b_{k-1}) & \text{if $i=k-1$,} \\
C(b_k,b_{k-1}) \Big( A_{k-1,i}(b_1,\ldots,b_{k-2},b_{k}) & \\
\qquad\qquad-A_{k-1,i}(b_1,\ldots,b_{k-1})\Big) & \text{if $1\leq i\leq k-2$.}
\end{cases}
\end{equation}
Finally combine~\eqref{rhs} with~\eqref{basic} and~\eqref{A k-i} to deduce the desired identity.
\section{Concluding remarks}\label{sec conclusion}
For simplicity of notation,
let $F(y,z)$ denote the left hand side of Corollary~\ref{cor k-2,3}(a), let
$K(x,y,z)$ denote the left hand-side of Corollary~\ref{cor k-2,3}(b),
and let $L(x,y,z)$ be the dual of $K(x,y,z)$ given by
\[
L(x,y,z) = \sum_{n,m\geq 0}\fr{q^{2n+m}(z,z^{-1};q)_{n} (y,y^{-1};q)_{n+m}}{(qy,q/y;q)_{n}(qx,q/x;q)_{n+m}}.
\]
Then it is natural to ask for a companion of
Corollary~\ref{cor k-2,3}(b) dealing with the $L(x,y,z)$.
We have the following connection between $K(x,y,z)$ and $L(x,y,z)$
\[
K(x,y,z)
=\sum_{m,n\geq 0}\fr{q^{2m+n} (y,y^{-1};q)_m (z,z^{-1};q)_{m+n}}{(qx,q/x;q)_m (qy,q/y;q)_{m+n}}
\]
\[
=\sum_{m\geq 0}\fr{q^m (y,y^{-1};q)_m}{(qx,q/x;q)_m}
\sum_{n\geq m}\fr{q^n(z,z^{-1};q)_{n}}{(qy,q/y;q)_{n}}
\]
\[
= \sum_{m\geq 0}\fr{q^m (y,y^{-1};q)_m}{(qx,q/x;q)_m}
\sum_{n\geq 0}\fr{q^n(z,z^{-1};q)_{n}}{(qy,q/y;q)_{n}}
\]
\[
-\sum_{m\geq 0}\fr{q^m (y,y^{-1};q)_m}{(qx,q/x;q)_m}
\sum_{n=0}^{m-1}\fr{q^n(z,z^{-1};q)_{n}}{(qy,q/y;q)_{n}},
\]
\[
= F(x,y) F(y,z)
-\sum_{m\geq 0}\fr{q^m (y,y^{-1};q)_m}{(qx,q/x;q)_m}
\sum_{n=0}^{m}\fr{q^n(z,z^{-1};q)_{n}}{(qy,q/y;q)_{n}}
+ \sum_{m\geq 0} \fr{q^{2m}(y,y^{-1};q)_m (z,z^{-1};q)_m}{(qx,q/x;q)_m (qy,q/y;q)_m}.
\]
\begin{equation}\label{K-L}
= F(x,y) F(y,z)-L(x,y,z)+ \sum_{m\geq 0} \fr{q^{2m}(y,y^{-1};q)_m (z,z^{-1};q)_m}{(qx,q/x;q)_m (qy,q/y;q)_m}.
\end{equation}

From~\eqref{K-L} the desired companion can be achieved by evaluating
the single sum
\[
 \sum_{m\geq 0} \fr{q^{2m}(y,y^{-1};q)_m (z,z^{-1};q)_m}{(qx,q/x;q)_m (qy,q/y;q)_m}
\]
as a linear combination of infinite products.

\bigskip
\noindent{\bf Acknowledgment.} The authors are  grateful to the referee for valuable comments and interesting suggestions which have improved the presentation and quality of the paper.

\noindent{\bf Data Availability Statement.\ }
Not applicable.
\end{document}